\documentclass[12pt]{article}

\usepackage[numbered,framed]{matlab-prettifier}
\usepackage{enumitem}
\usepackage{listings}
\usepackage{adjustbox}
\usepackage[margin=1in,footskip=0.25in]{geometry}
\usepackage{relsize}
\usepackage{amsthm}%
\usepackage{mathtools}
\usepackage{amsfonts}%
\usepackage{pdfpages}
\usepackage{amsmath}
\usepackage{color}
\usepackage{flexisym}

\definecolor{mygreen}{RGB}{28,172,0} 
\definecolor{mylilas}{RGB}{170,55,241}

\usepackage{cite}
\newtheorem{theorem}{Theorem}
\newtheorem{remark}[theorem]{Remark}

\begin{document}

\title{Fractional order logistic map: Numerical approach}

\author{Marius-F. Danca{\footnote{Corresponding author}}\\
Romanian Institute of Science and Technology, \\
400487 Cluj-Napoca, Romania\\
Email: danca@rist.ro\\
}

\maketitle

\begin{abstract}
In this paper the fractional order logistic map in the sense of Caputo's fractional differences is numerically approached. It is shown that the necessary iterations number to avoid transients must be of order of thousand, not of order of hundreds as commonly used in several works. Also, it is revealed an interesting phenomenon according to which for every initial condition it correspond a different bifurcation diagram. This phenomenon seems to appear also in other Fractional Order (FO) difference systems, fact which could represent an obstacle for the numerical analysis. A short Matlab code is used to obtain the results.

\end{abstract}

\textbf{keyword }Caputo fractional differences; FO logistic map; Bifurcation diagram

\vspace{3mm}

\section{Introduction}

The concept of derivatives of non-integer values dates back to the beginning of the theory of differential calculus, and the development of the theory of fractional calculus dates from the work of Euler, Liouville, Riemann, Letnikov \cite{01,02}. The fractional derivatives and integrals are useful in engineering and mathematics, being helpful for scientists and researchers working with real-life applications (see, e.g., \cite{03}). For basic aspects of the theory of fractional differential equations refer the monograph \cite{04}, or \cite{00} for a review of definitions of fractional derivatives and other operators.

On the other hand, fractional difference equations received increasing attention recently, due the description of memory and hereditary properties, and one of the first definitions of a fractional difference operator has been proposed in 1974 \cite{04}. Problems related to the left and right Caputo fractional sums and differences are presented in \cite{patru}. However, there still are only few works in the theory of the fractional finite difference equations \cite{frac20}.
Initial value problems are studied in \cite{ati}.
The stability problem is studied in \cite{unu,doi,trei}.
Some applications have been provided in \cite{cinci,sase,sapte}

An appropriate bibliography for the fractional $q$-calculus can be found in \cite{ati2}.
Some discrete fractional maps, such as logistic map, variant of logistic map and standard map are presented in \cite{opt,noua,zece,unspe,doispe,treispe}. In \cite{patrupatru} the chaos, 0-1 test, the $C_0$ complexity, entropy, and the control of the discrete fractional Duffing system are studied. In \cite{treitrei} the convergence speed depending on the initial conditions, is analyzed. Synchronization of fractional-order discrete-time chaotic systems and application to secure communication are presented in \cite{cincicinci} and \cite{sasesase}. Fractional-order discrete-time uncertain systems and their LQ optimal control proposed in \cite{saptesapte}. In \cite{optopt} the stability of nonlinear discrete fractional systems with application to vibrating eardrum.

In this paper we extend the range of iterations of the $q$-order left fractional differences of the logistic map with $q\in(0,1)$, in order to underline some new phenomena appearing in bifurcation diagrams versus the fractional order and the bifurcation parameter. It is shown that only hundreds iterations are not enough to characterize numerically the dynamics of the fractional logistic map. First, the form of the initial value problem of the fractional order logistic map and its numerical integration is presented. Next, the numerical approach is realized, mainly, via bifurcation diagrams versus the fractional order and versus the bifurcation parameter.

\section{The FO logistic map}

To obtain the fractional order form of the logistic map (FOLM) in Caputo's left fractional (with delta operator $\Delta$) differences sense \cite{patru}, let $N_a=\{a,a+1,a+2,\ldots\}$. Then, for~$q>0$ and $q\not\in \mathbb{N}$, the~$q$-th Caputo-like discrete fractional difference with starting point $a$ of a function $u:N_a\rightarrow \mathbb{R}$ is defined as~\cite{ati,xxx2}
\begin{equation*}\label{capa}
\Delta_a^q u(t)=\Delta_a^{-(m-q)}\Delta^m u(t)=\frac{1}{\Gamma(m-q)}\sum_{s=a}^{t-(m-q)}(t-s-1)^{(m-q-1)}\Delta^mu(s),
\end{equation*}
for $t\in N_{a+m-q}$ and $m=[q]+1$.

$\Delta^m$ is the $m$-th order forward difference operator,
\[
\Delta^m u(s)=\sum_{k=0}^{m}\binom {n}{k}(-1)^{m-k}u(s+k),
\]
while $\Delta_a^{-q}$ represents the $q$-th fractional sum of $u$ starting from $a$, namely,
\begin{equation}\label{suma}
\Delta_a^{-q}u(t)=\frac{1}{\Gamma(q)}\sum_{s=a}^{t-q}(t-s-1)^{(q-1)}u(s),~t\in \mathbb{N}_{a+q},
\end{equation}
with the falling factorial $t^{(q)}$ in the following form:
\[
t^{(q)}=\frac{\Gamma(t+1)}{\Gamma(t-q+1)}=t(t-1)...(t-q+1).
\]

Note that the fractional operator $\Delta_a^{-q}$ maps functions defined on $\mathbb{N}_a$ to functions on~$\mathbb{N}_{a+q}$, i.e. in \eqref{suma} $u$ is defined for $t=a~\mathrm{mod}(1)$, and $\Delta_a^{-q} u$ is defined for $(a+q)~\mathrm{mod}~ (1)$.

For the case considered in this paper, $q\in(0,1)$, when $m=1$, $\Delta u(s)=u(s+1)-u(s)$, and~Caputo's fractional difference, denoted hereafter $\Delta_*^q$, becomes
\begin{equation*}
\Delta_*^q u(t)=\frac{1}{\Gamma(1-q)}\sum_{s=a}^{t-(1-q)}(t-s-1)^{(-q)}\Delta u(s).
\end{equation*}

Considering the usual case of zero starting point in the fractional sum \eqref{suma}, $a=0$, one can consider the following Caputo's like discrete Initial Value Problem (IVP):
\begin{equation}\label{eq2}
{\Delta_*^q u(t)=f(t+q-1),u(t+q-1)),~ t\in N_{1-q}, ~u(0)=u_0,
}\end{equation}

\noindent {for $q\in(0,1)$ and $f$ is a real-valued function and $u_0$ a real number.}

The solution of the IVP \eqref{eq2} is given by~\cite{xxx2}
\begin{equation}\label{inte}
u(t)=u_0+\frac{1}{\Gamma(q)}\sum_{s=1-q}^{t-q}(t-s-1)^{(q-1)}f(u(s+q-1)).
\end{equation}

Note that the recursive iteration in the sum equation \eqref{inte} implies that the solution is unique.

To obtain a convenable numerical form, consider in \eqref{inte} the following substitution $s+q=i$. Then, the falling factorial $(t-s-1)^{(q-1)}$ becomes
\[
(t-s-1)^{(q-1)}=\frac{\Gamma(t-s)}{\Gamma(t-s-q)}=\frac{\Gamma(t-1+q)}{\Gamma(t-s-q+1)}=\frac{\Gamma(t-i+q)}{\Gamma(t-i+1)},
\]
and, because $t\in\mathbb{N}_1=\{1,2,...\}$, by replacing $t\in \mathbb{N}$ with the usual index $n\in \mathbb{N}$, a convenient iterative numerical form of the integral \eqref{inte}, which will be used in this paper, is
\begin{equation}\label{eq3}
u(n)=u(0)+\frac{1}{\Gamma(q)}\sum_{i=1}^n\frac{\Gamma(n-i+q)}{\Gamma(n-i+1)}f(u(i-1)),~ u(0)=u_0, ~{n\in \mathbb{N}_1~}.
\end{equation}

Now, consider the IVP of the FOLM

\begin{equation}\label{eqx}
{\Delta_*^q x(t)=px(t+q-1)(1-x(t+q-1)),~ t\in N_{1-q}, ~x(0)=x_0,
}\end{equation}
with $q\in(0,1)$ and where the ranges of $x$ and $p$ will no longer be as for the integer order case.

The numerical solution \eqref{eq3}, becomes
\begin{equation}\label{eqq}
x(n)=x(0)+p\frac{1}{\Gamma(q)}\sum_{i=1}^n\frac{\Gamma(n-i+q)}{\Gamma(n-i+1)}x(i-1)(1-x(i-1)),~ x(0)=x_0, ~{n\in \mathbb{N}_1~}.
\end{equation}

Note that the integral \eqref{eqq} being unique, defines a discrete dynamical system \cite{xxx} which fully determines the behavior of the FOLM and therefore we are entitled to study \eqref{eqq} in order to characterize the dynamics of the FOLM.

\section{Numerical approach of the FOLM}

In this section the FOLM is numerically analyzed on Bifurcation Diagrams (BDs), via the numerical scheme \eqref{eqq} The initial conditions are chosen empirically. To perform the numerical analysis it is necessary to increase the iterations number, $n_{max}$ which, usually in the existing works, are of order of hundreds only and do not represent the real behavior of the FOLM.

In this paper the iterative solution \eqref{eqq} is implemented in the Matlab and, therefore, it is necessary to avoid the annoying zero index problem in \eqref{eqq}. One of the several easy, but time consuming ways, is the code \cite{cod} which allows to deal with zero-indexable arrays. In this paper, a simple code with a trick to overcome the zero index, is presented in Appendix. Note that the code can be used for other functions. Also, the speed of the code can be significantly improved using different Matlab implemented functions.

\begin{figure}
\begin{center}
\includegraphics[scale=0.65]{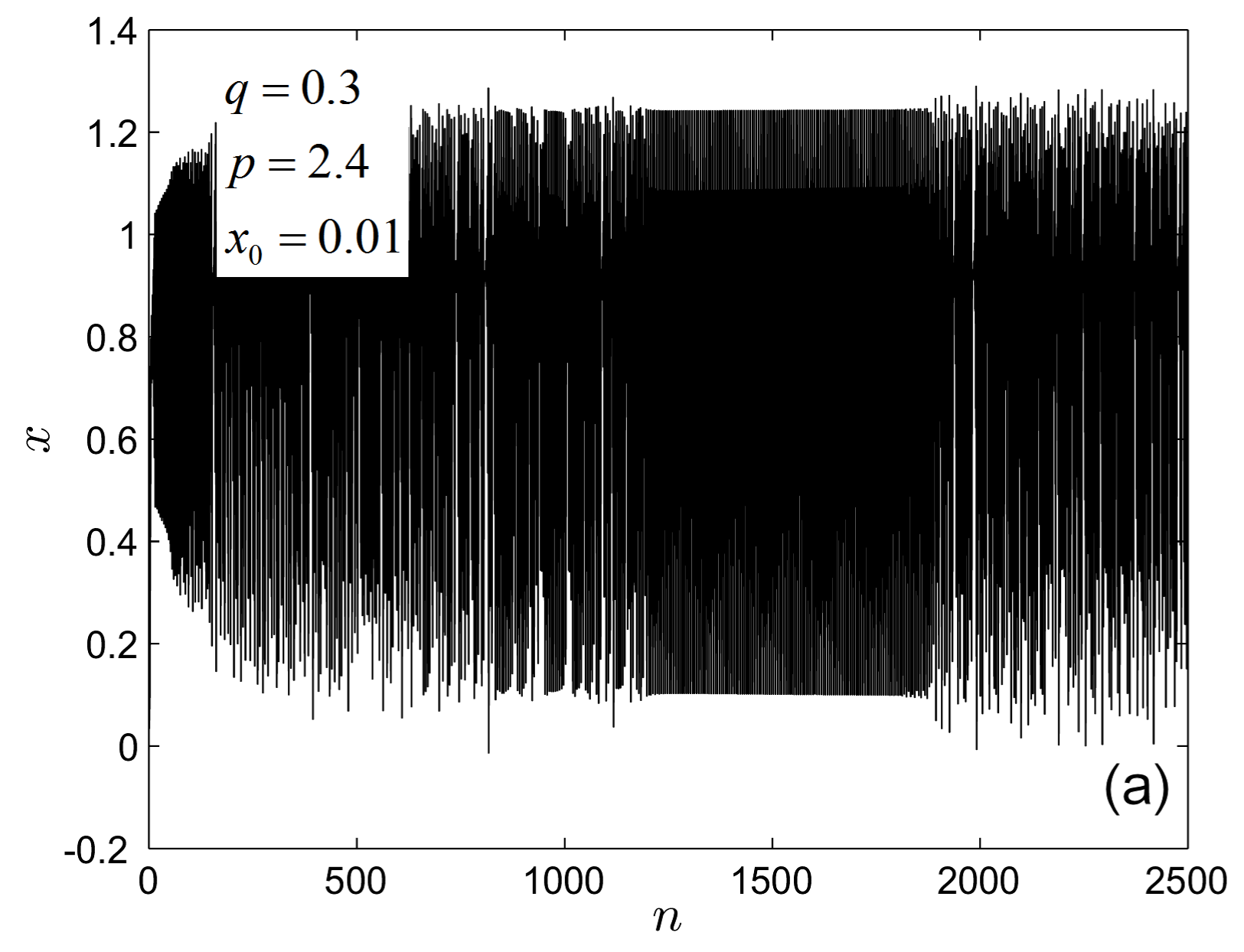}
\caption{Time series of the FOLM for $q=0.3$ and $p=2.4$. For $n\in[0,N]$ with $N=2500$, two sets of transients (chaotic and periodic like) appear: chaotic for $n\in[0,1100)$ and periodic-like for $n\in(1100,1800)$.}
\label{fig000}
\end{center}
\end{figure}

The BDs of the FOLM are constructed for several thousands of iterations $n_{max}$: 2500, 5000 and 7500. Fewer iterations, of order of hundreds, could lead to wrong solutions by containing only transients.
For example in Fig. \ref{fig000} the time series for $q=0.3$ and $p=2.4$ is presented. For values of $n_{max}$ less than about 1100, the system appears to be chaotic, while for $n_{max}\in(1100,1800)$ the behavior appears periodic. However, these dynamics are only transitory and, for this case, the real behavior can be considered only for $n_{max}>1800$. Therefore, in this paper the choice $n_{max}\geq2500$ has been considered a reasonable compromise between neglecting the transients and computational time.

In order to increase the value of $n_{max}$, in the solution \eqref{eqq} one can uses the following identity which allows to extend $n_{max}$ from few hundreds to several thousands
\[
\frac{\Gamma(n-i+q)}{\Gamma(n-i+1)}=e^{\ln\Gamma(n-i+q)-\ln\Gamma(n-i+1)}.
\]

\begin{figure}
\begin{center}
\includegraphics[scale=0.65]{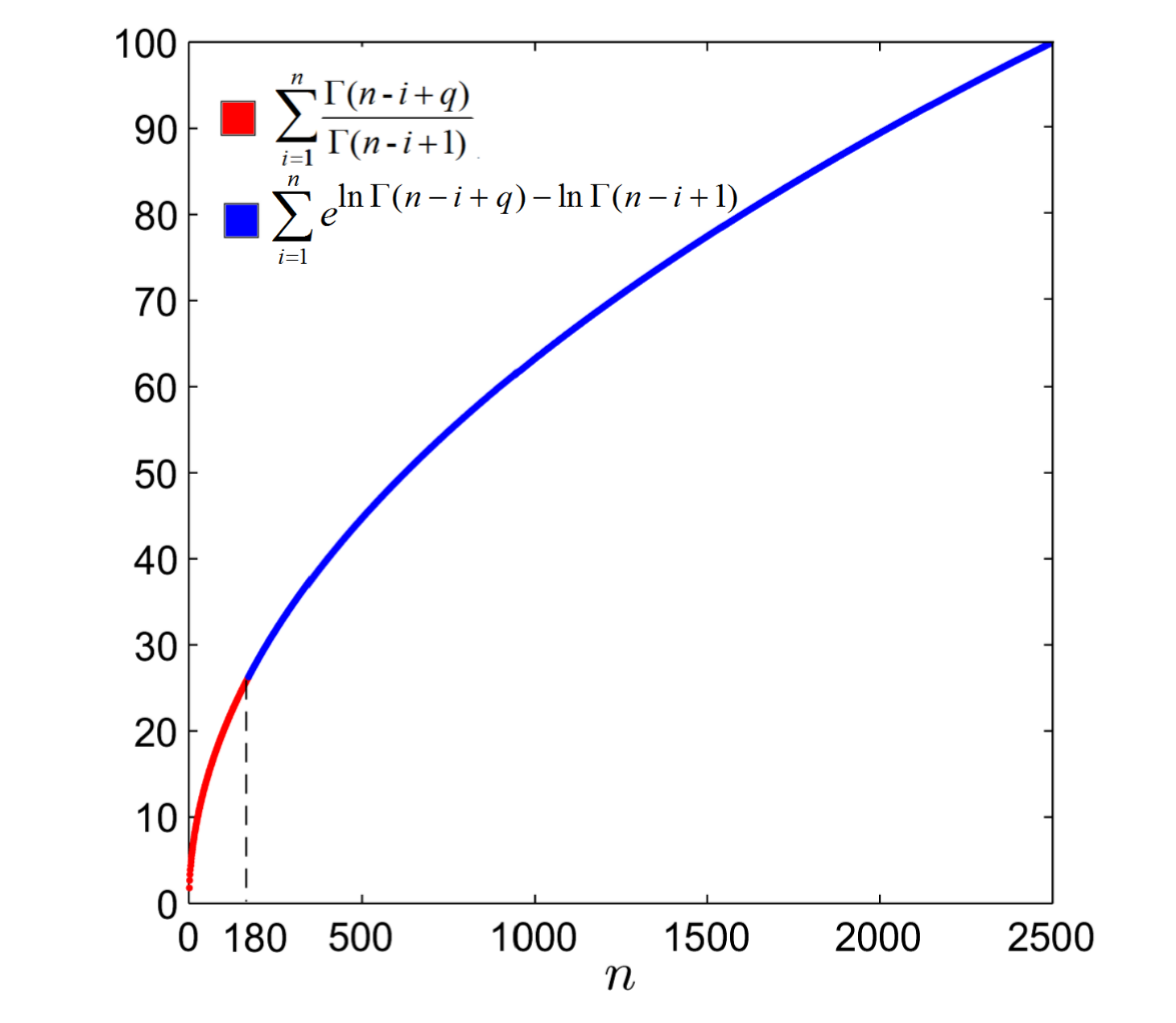}
\caption{Graphs of the sums $\sum_{i=1}^n\frac{\Gamma(n-i+q)}{\Gamma(n-i+1)}$ (red plot) and $\sum_{i=1}^n\exp^{ln \Gamma(n-i+q)-ln\Gamma(n-i+1)}$ (blue plot).}
\label{fig00}
\end{center}
\end{figure}

The difference between the two approaches is evident in Fig. \ref{fig00}, where the sum in \eqref{eq3} is presented. Because $f$ is bounded, the sum can be calculated with $f$ constant (e.g. 1), in the two ways (red plot with $\frac{\Gamma(n-i+q)}{\Gamma(n-i+1)}$ and blue plot with $e^{\ln\Gamma(n-i+q)-\ln\Gamma(n-i+1)}$ ). As can be seen, for $n_{max}>180$, $\frac{\Gamma(n-i+q)}{\Gamma(n-i+1)}$ becomes unbounded.

\begin{figure}
\begin{center}
\includegraphics[scale=0.5]{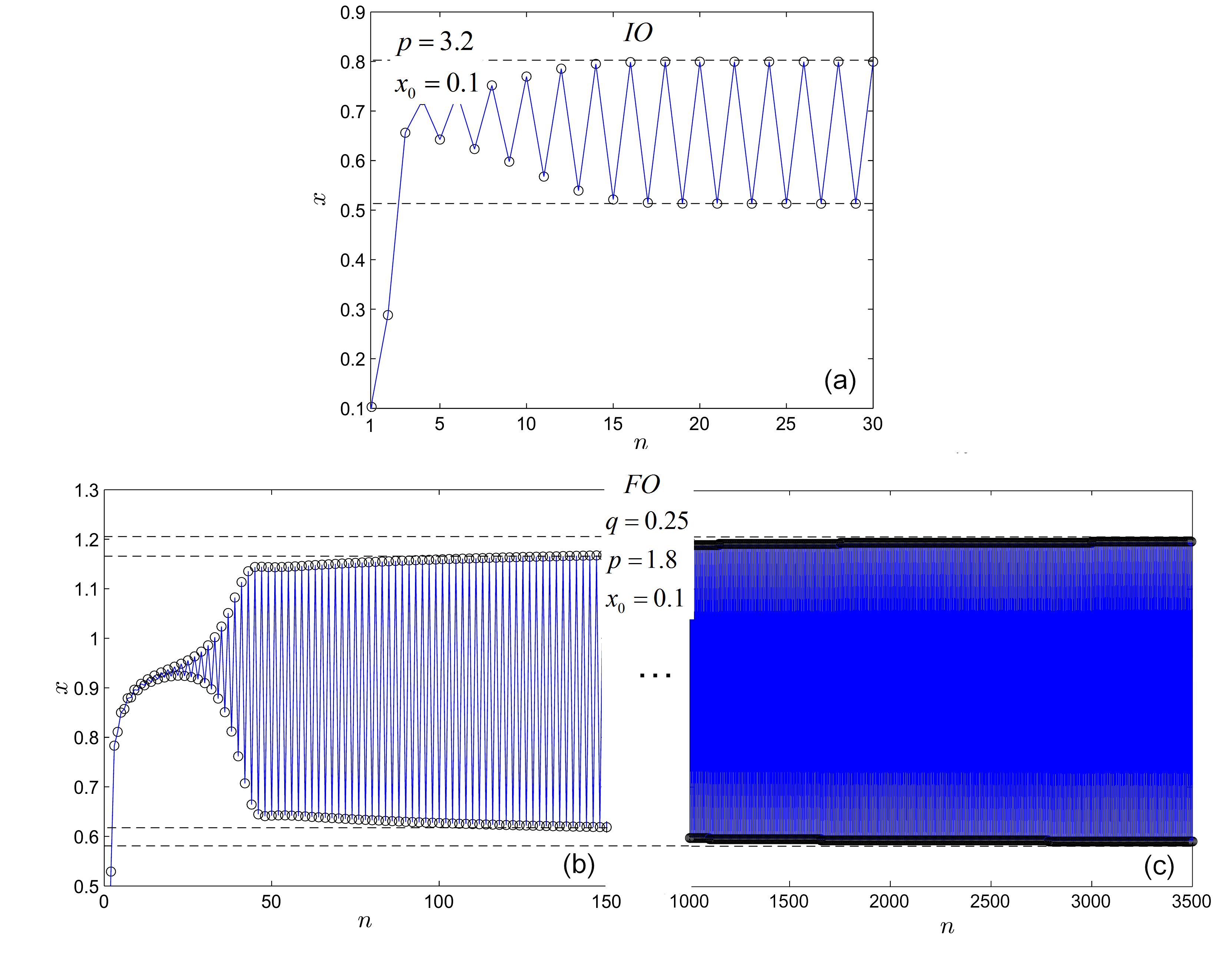}
\caption{Time series of IOLM and FOLM. (a) Time series revealing the $2$-period stable cycle of the IOLM for $p=3.2$ and $x_0=0.1$. A short preperiod is necessary only to reach the cycle; (b)-(c) Time series of the FOLM for $q=0.25$, $p=1.8$ and $x_0=0.1$. The $2$-NPO reaches asymptotically a periodic orbit. }
\label{fig0}
\end{center}
\end{figure}

\begin{remark}\label{remus}
Since FO systems (continuous or discrete) admit not exactly periodic solutions, the orbits of the FOLM \eqref{eqx} which ``look'' periodic, are only apparent periodic (see e.g. \cite{neper1} for continuous systems and \cite{nepe} for discrete systems). Beside the existing analytical proofs of this property, the long-time history typical to all kind of FO systems underlines logically the non-periodicity of these systems. For example, in Figs. \ref{fig0} are presented comparatively the case of Integer Order Logistic Map (IOLM) $x(n+1)=px(n)(1-x(n))$, $n\in \mathbb{N}$, with $p=3.2$ for which, after a short transients (preperiod), the system behaves along a 2-period stable cycle (Fig. \ref{fig0} (a)). In the case of the FOLM, for $q=0.25$ and $p=1.8$ (Figs. \ref{fig0} (b), (c)), it can be seen that the periodicity is apparently even after $n_{max}=3500$ iterations and the transients last infinitely. In this paper, these kind of orbits are called numerically-periodic orbits (NPOs) \cite{000,neper2}.
\end{remark}

\begin{figure}
\begin{center}
\includegraphics[scale=1]{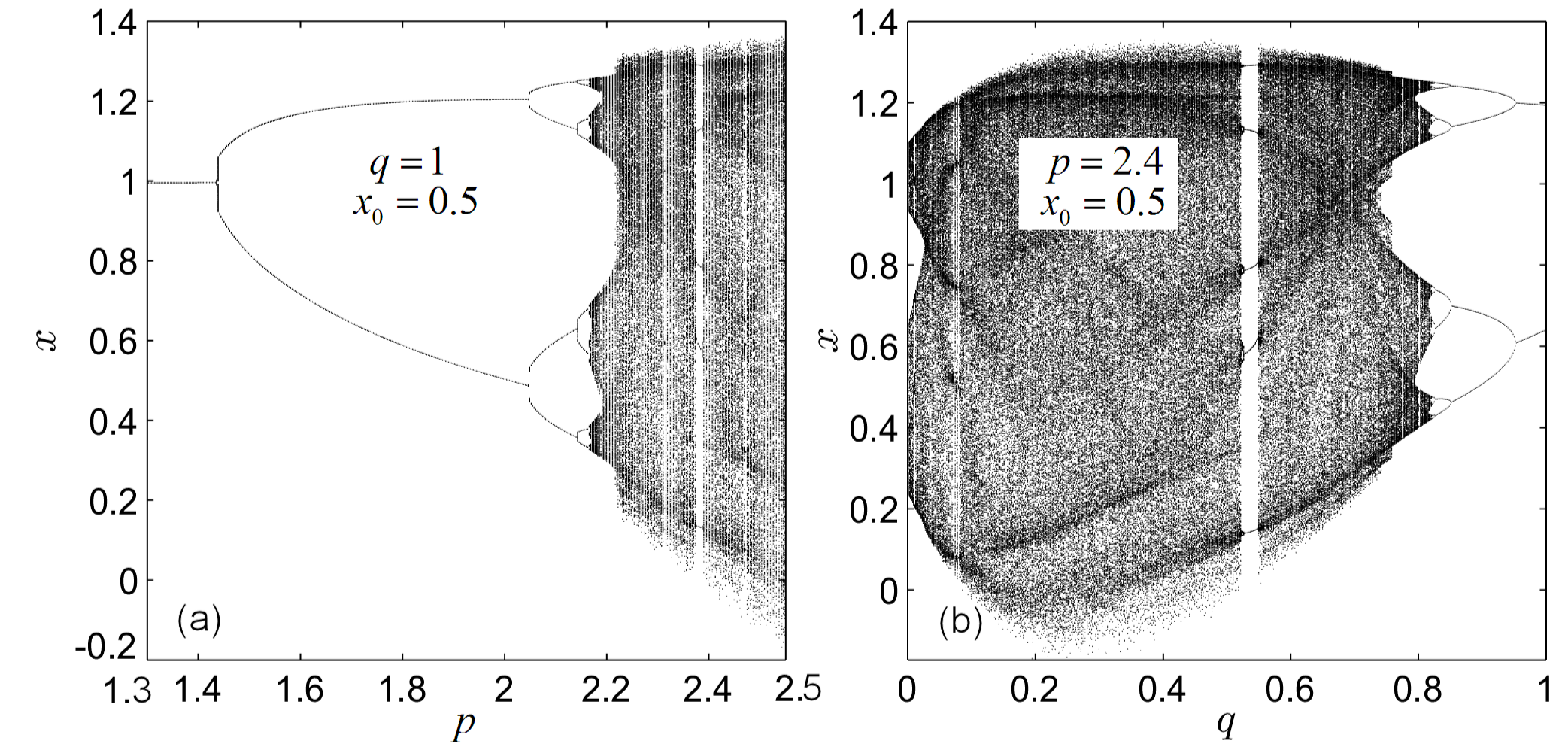}
\caption{BDs of the FOLM. (a) BD vs. $p\in[1.3,2.5]$ for $q=1$; (b) BD vs. $q\in(0,1)$ for $p=2.4$.}
\label{fig1}
\end{center}
\end{figure}

Consider for $x_0=0.5$ the BDs of the FOLM vs $p\in[1.3,2.5]$ and $q=1$ (Fig. \ref{fig1} (a)) and vs $q\in(0,1)$ for $p=2.4$ (Fig. \ref{fig1} (b)).

\begin{figure}
\begin{center}
\includegraphics[scale=0.6]{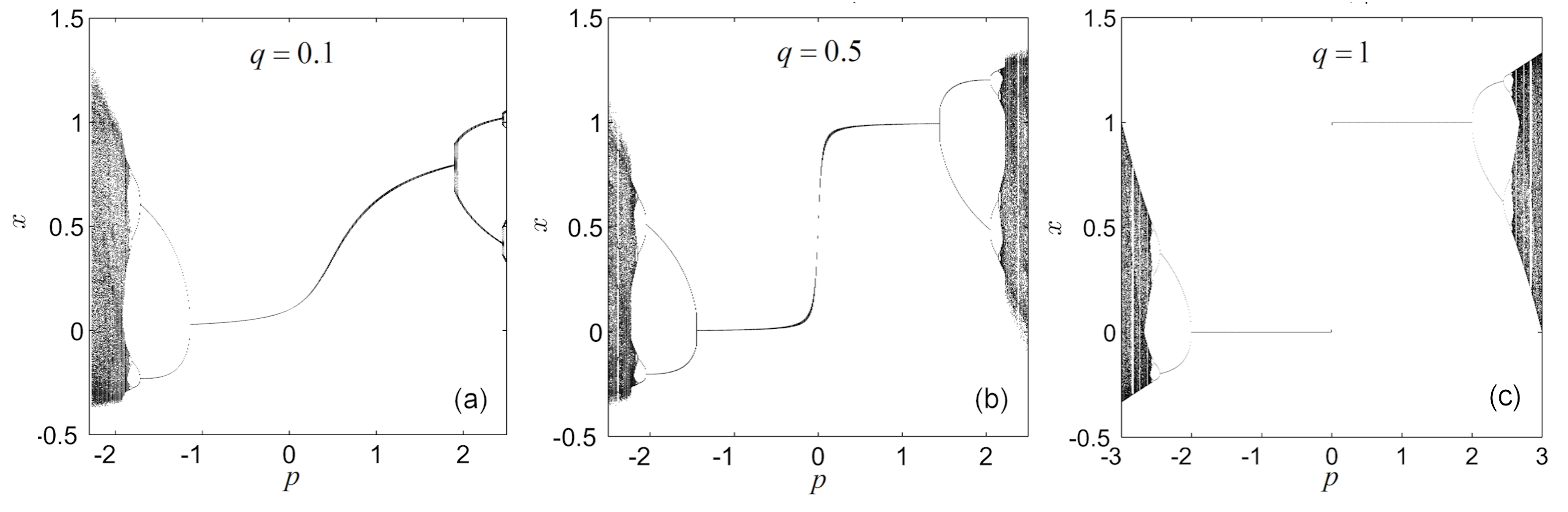}
\caption{BDs of the FOLM vs. $p$ for a single initial condition and different values of $q$. (a) Case $q=0.1$; (b) Case $q=0.5$; (c) Case $q=1$.}
\label{fig2}
\end{center}
\end{figure}

Even the shape of the BD vs $p$ in Fig. \ref{fig1} (a) is a reminiscent of the case of IO, there are significant differences such as the values of the state variable $x$. Also, while in the IO case $p$ is usually chosen within the interval $[0,4]$, in the FO case, the admissible range of $p$ could be significantly larger with spectacular behaviors as shown in Figs. \ref{fig2} where $q=0.1$ vs $p\in[-2.5,2.5]$ (Fig. \ref{fig2} (a)), $q=0.5$ (Fig. \ref{fig2} (b)) and $q=1$ vs $p\in[-3,3]$ (Fig. \ref{fig2} (c)). In Figs. \ref{fig2} (a), (b), the maximal parameter interval is $[-2.5,2.5]$, while for $q=1$ (Fig. \ref{fig2} (c)) the interval could be larger $[-3,3]$. Beside the symmetry which begins at $q=0.5$ and continues till $q\rightarrow 1$ (Figs. (b), (c)), it can be seen that for smaller values of $q$ ($q<0.5$, Fig. \ref{fig2} (a)), chaos exists only for negative values of $p$.

However, the most important differences between IO and FO logistic map can be unveiled by considering BDs as generated from several initial conditions whose generated patterns tending to fixed points, periodic orbits, or chaotic attractors (in the sense of Remark \ref{remus}) are overplotted in the space $(p,x)$, or $(q,x)$, in different colors (technique used in finding attractors coexistence).

For a better understanding of underlying phenomenons in these cases, denote the bifurcation set corresponding to one of the considered initial conditions by Bifurcative Set (BS). To each initial condition corresponds a BS and, therefore, a BD made with, e.g., three or more initial conditions, is composed by three or more BSs which could look ``similar'' or different.
\emph{Similarity }is considered in this paper as an identity-like of the shapes of BSs in the spaces $(p,x)$ or $(q,x)$.

\begin{figure}
\begin{center}
\includegraphics[scale=0.9]{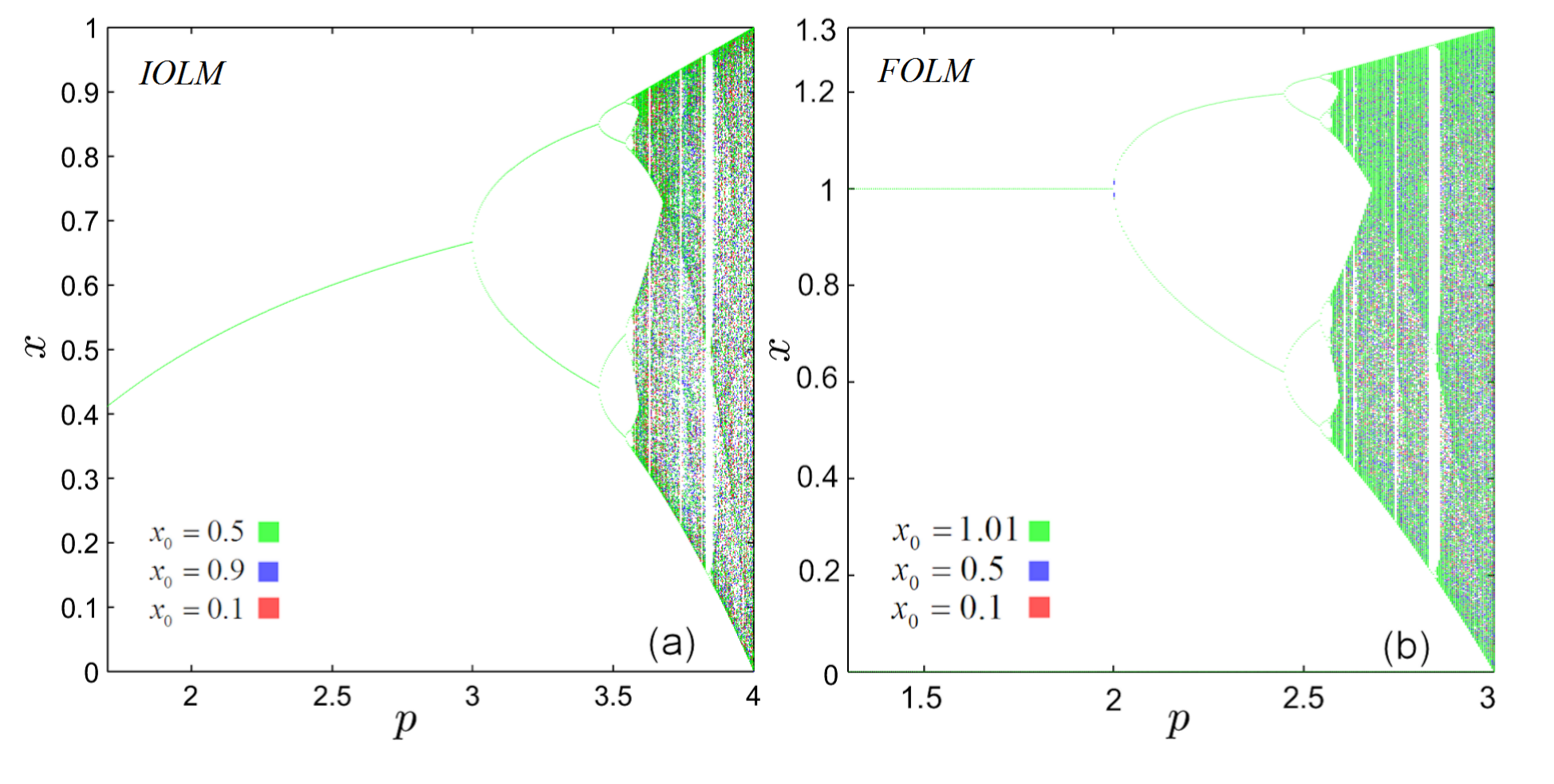}
\caption{BDs of the IOLM and FOLM. (a) BD of the IOLM composed by three BSs generated with three initial conditions $x0=0.5$, $x_0=0.9$ and $x_0=0.1$ (green, blue and red, respectively); (b) BD of the FOLM composed by three BSs generated with three initial conditions $x0=1.01$, $x_0=0.5$ and $x_0=0.1$ (green, blue and red, respectively). }
\label{fig3}
\end{center}
\end{figure}

As expected and well known, all initial conditions in the case of IOLM, generate similar BSs. In Fig. \ref{fig3} (a) three initial conditions, $x_0=0.5$, $x_0=0.9$ and $x_0=0.1$, are considered, the underlying BSs being colored green, blue and red respectively. Similarly, for the limit case $q=1$ (note that \eqref{eqq} is deduced for $q\in(0,1)$), the FOLM with initial conditions, $x_0=1.01$, $x_0=0.5$ and $x_0=0.1$, presents similar BSs. However, this similarity stops at the case $q=1$ and for $q<1$ unexpected behaviors appear.

\begin{figure}
\begin{center}
\includegraphics[scale=1]{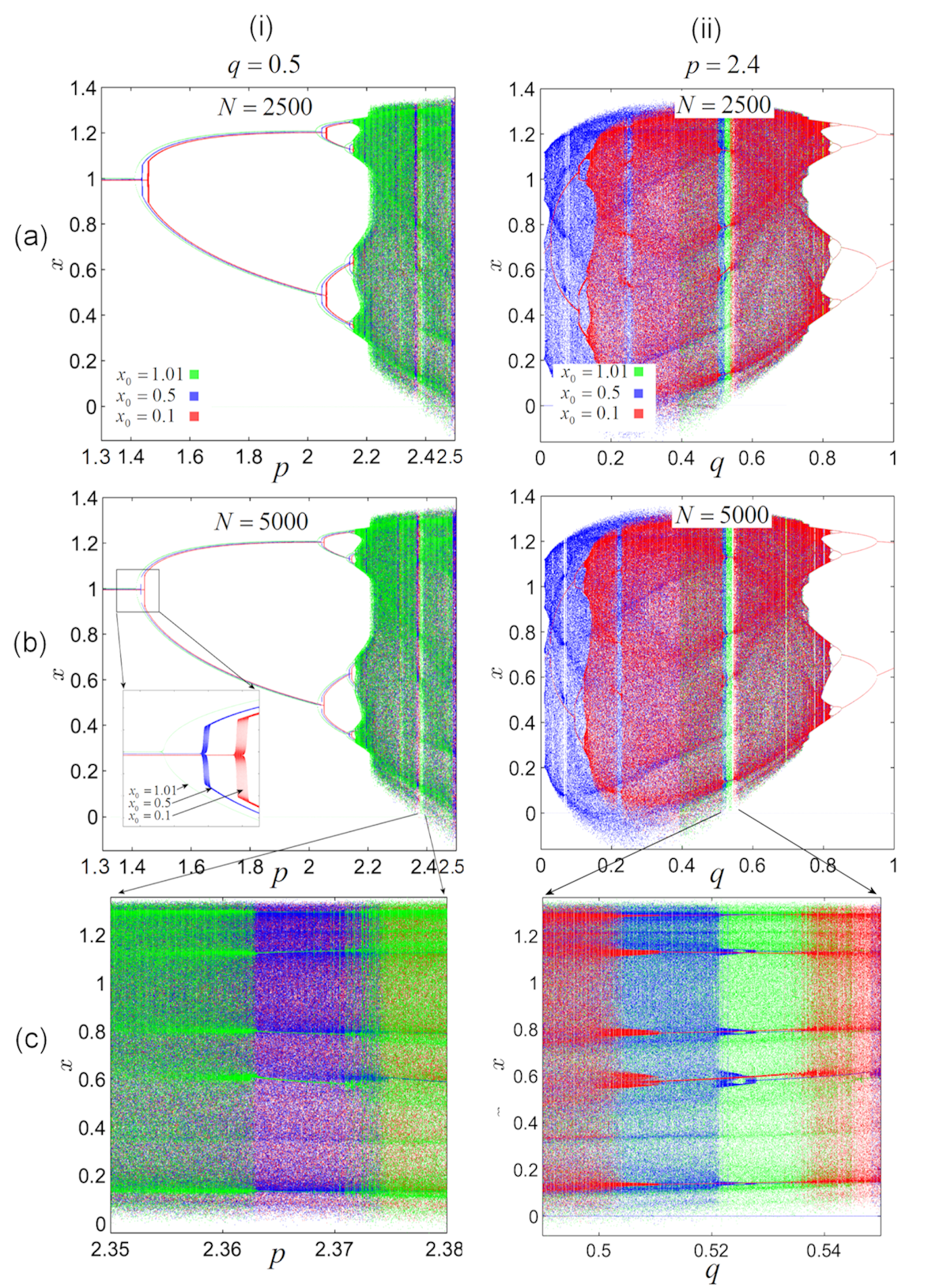}
\caption{BDs of the FOLM. Column (i): BD vs. $p$ for $q=0.5$ composed by three BSs for initial conditions $x0=1.01$, $x_0=0.5$ and $x_0=0.1$ (green, blue and red, respectively); Column (ii): BD vs. $q$ for $p=2.4$ composed by three BSs generated from the same initial conditions; Line (a): $N_{max}=2500$; Line (b): $n_{max}=5000$; Line (c): Zooms of line (b).}
\label{fig4}
\end{center}
\end{figure}

Consider the three initial conditions $x_0=1.01$, $x_0=0.5$ and $x_0=0.1$ and the BSs for $q=0.5$ and $p\in[1.3,2.5]$ (column (i) in Figs. \ref{fig4}) and the BSs for $p=2.4$ and $q\in(0,1)$ (column (ii) in Figs. \ref{fig4}). Note that in the space $(q,x)$, the BS obtained from the initial condition $x_0=1.01$ exists only for about $q>0.4$, the orbits for $q\in(0,0.4)$ being divergent.

As can be seen, all BSs are different one to each other along axes $p$ and $q$. Thus, for $n_{max}=2500$ iterations, one can see in the BD vs $p$ (Fig. \ref{fig4} (a), column (i)) that the first bifurcation point multiplies (to the three initial conditions, correspond three bifurcation points). This translation-like of the bifurcation points in the space $(p,x)$ repeats at every bifurcation points. Also, the chaotic bands in each BSs differ. In Fig. \ref{fig4} (a), column (ii), one can see that the differences between BSs are even more significant. While for some ranges of $q$, such as $q=0.1$, the BS corresponding to $x_0=0.5$ (blue plot) indicates the existence of chaos, the BS corresponding to $x_0=0.1$ suggests the existence of NPOs (red plot). Also, an interesting situation is unveiled: for $q$ increasing to $1$, the mentioned differences vanish, fact which is consistent with the case $q=1$ presented in Fig. \ref{fig3} (a).

\begin{figure}
\begin{center}
\includegraphics[scale=0.95]{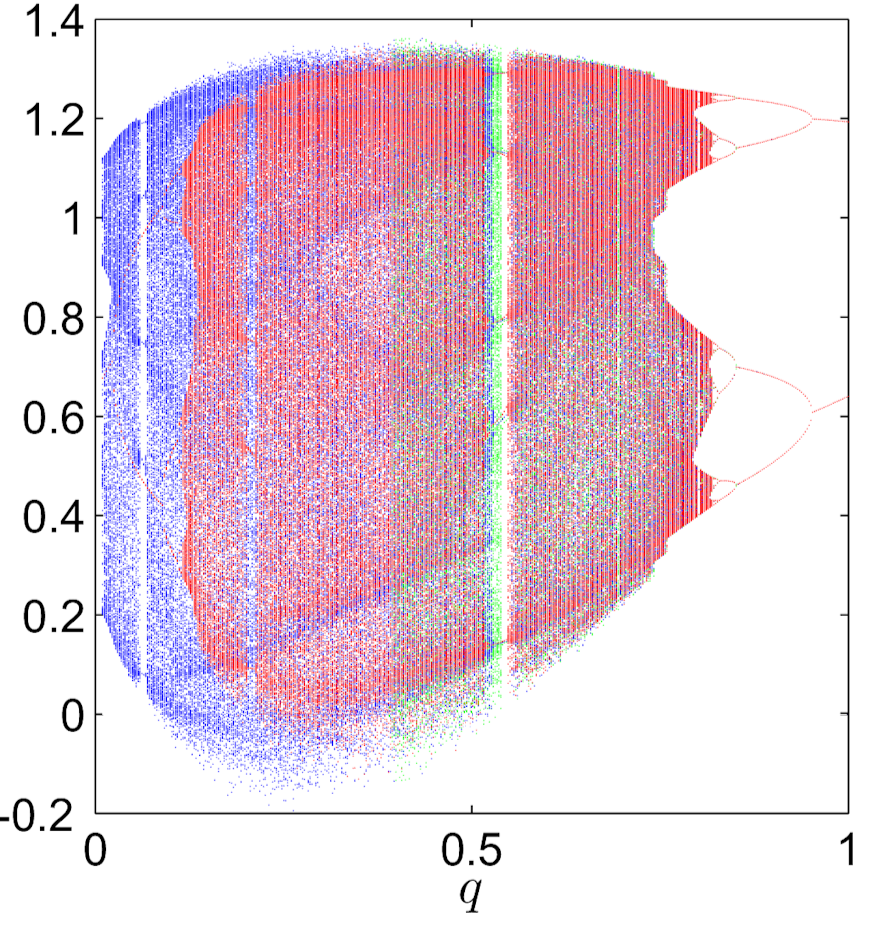}
\caption{BD of the FOLM vs. $q$ for $a=2.4$ with $n_{max}=7500$.}
\label{fig44}
\end{center}
\end{figure}

At this point a natural question arises: are $n_{max}=2500$ iterations sufficient enough to unveil the ``real'' shapes of a BSs? In order to obtain an eloquent answer to this question, $5000$ iterations have been considered for both BDs, vs. $p$ and $q$ in Figs. \ref{fig4} (b), columns (i) and (ii). As can be seen, no significant differences between $n_{max}=2500$ and $n_{max}=5000$ can be remarked. Moreover, for larger values of $n_{max}$, the differences, smaller, still remain (see Fig. \ref{fig44} where the BD is made with $n_{max}=7500$ iterations for each of the BSs).

One can state numerically that the BSs of the FOLM depend on the initial conditions: regardless of the iterations number $n_{max}$, to every initial conditions in the spaces $(p,x)$ or $(q,x)$, correspond different BSs, $x_0\mapsto BS(x_0)$.

\begin{figure}
\begin{center}
\includegraphics[scale=0.7]{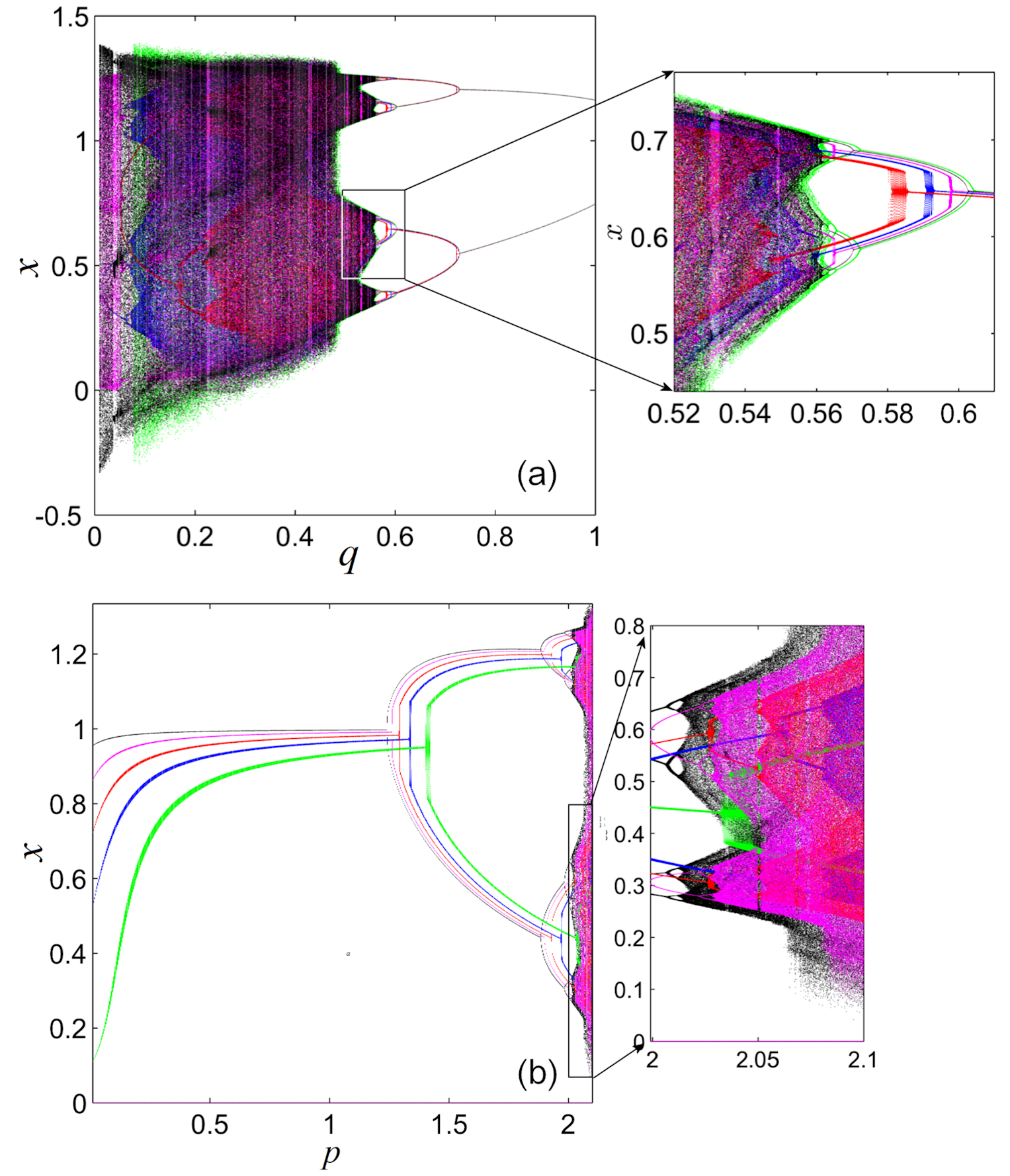}
\caption{BD of the FOLM composed by five initial conditions $(0.1),(0.5),(0.95),(0.7),(0.85)$. (a) BD vs. $q$ with $p=2.2$ and a zoomed region; (b) BD vs. $p$ with $q=0.3$ and a zoomed region. }
\label{fig5}
\end{center}
\end{figure}

In order to see if these phenomena appear for other parameters choice, the BDs case of 5 initial conditions $x_0\in\{0.1,0.5,0.95,0.7,0.85\}$ for $p=2.2$  are presented in Fig. \ref{fig5} (a) and for $q=0.3$ in Fig. \ref{fig5} (b). The zooms underline similar differences.

\begin{figure}
\begin{center}
\includegraphics[scale=0.75]{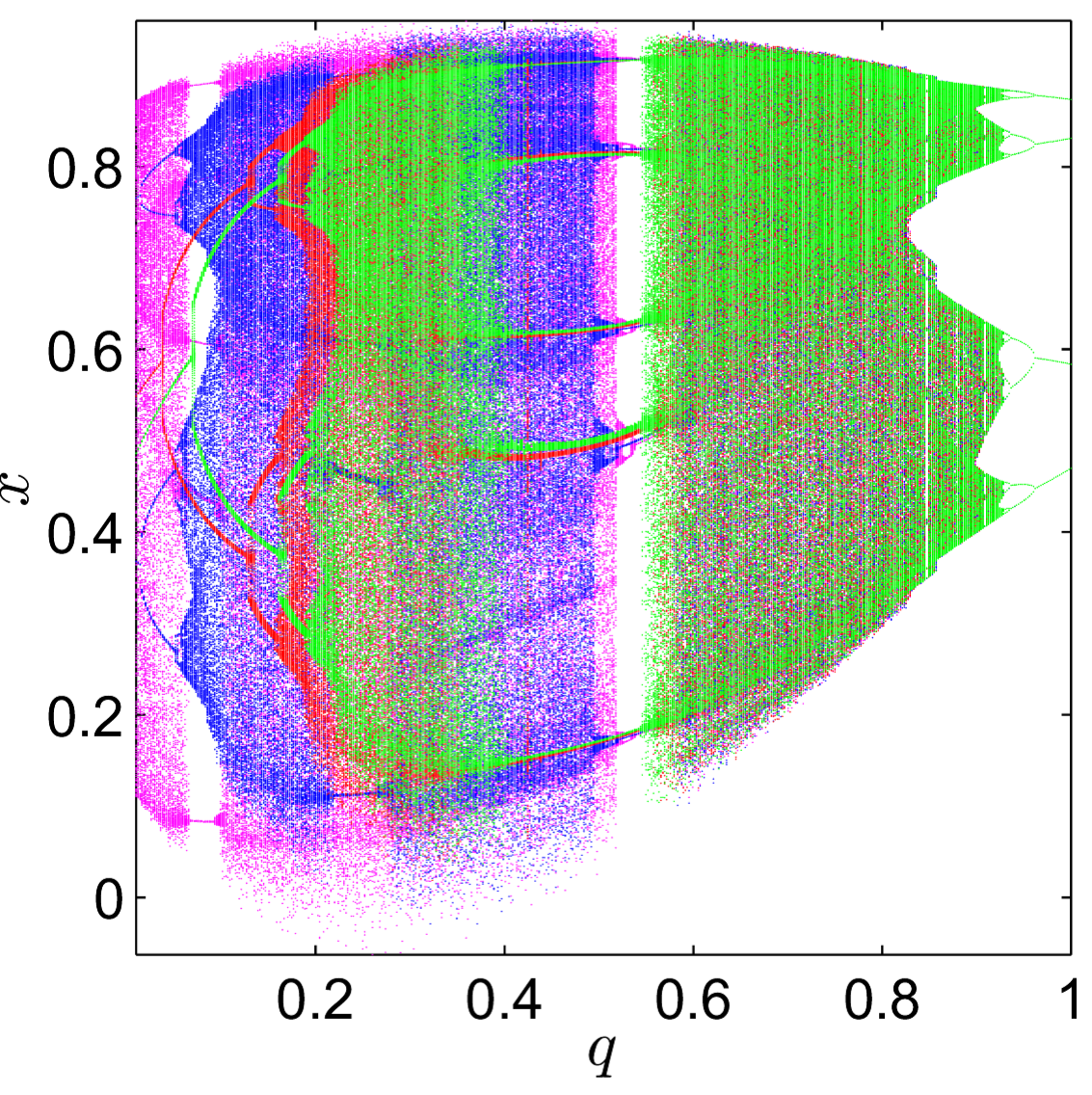}
\caption{The BD composed by four BSs in the case of FO Puu's system \eqref{pul} for $a=1.27$ and $x_0\in\{0.2,0.5,0.1,0.4\}$ (green, blue, red and magenta plot, respectively). }
\label{fig6}
\end{center}
\end{figure}

To note that other Caputo-like discrete systems of FO modeled by the IVP \eqref{eq2} seem to present similar phenomena. For example, in Fig. \ref{fig6} the BD composed by four BSs of the PUU system of FO \cite{dancus2} for $a=1.27$ and $x_0\in\{0.2,0.5,0.1,0.4\}$ is presented (note that the BD is symmetric about the $q$ axis and, therefore, only the upper part is considered here). The IO system is modeled by the recursive equation

\[
x(n + 1) = ax(n) -(a + 1)x^3(n), ~a\in \mathbb{R}_+,~ n\in \mathbb{N},
\]
and the IVO \eqref{eq2} is
\begin{equation}\label{pul}
\Delta^q_*x(n) = ax(n + q -1) -(a + 1)x^3 (n + q -1),~ n\in N_{1-q},~ x(0) = x_0.
\end{equation}

Beside the differences between BSs, the convergence-like of the BSs once $q$ approaches 1 can be seen, similarly with the case of the FOLM.

Summarizing, via the solution \eqref{eqq}, to every initial condition correspond different BSs. Thus, while for the IO case there exists an invariance of the BD with respect on initial conditions (including attractors coexistence), the FOLM is characterized by the non-invariance of the BD with respect the initial conditions.

Therefore, in the case of the FOLM, and possibly in other FO discrete systems modeled by the IVP \eqref{eq2}, the following issues arise:

\begin{itemize}[noitemsep,topsep=-2pt]
\item [i)] The cause of the differences could be a consequence of the mentioned long-time history.
\item[ii)] The computer errors cannot be considered as the basis of this phenomenon, since it appears only for about $q<0.5$, not for $q\in(0.5,1)$.
\item[iii)] How the BD of the FOLM can be obtained?
\item[iv)] As the zooms in Figs. \ref{fig4} (c), columns (i) and (ii) show, there exist windows revealing potential attractors coexistence. As known, an ingredient in finding numerically hidden attractors \cite{hid1} is the numerical analysis of their attraction basins to see if they intersect or not with small neighborhoods of equilibria (if the system admits equilibria, as in the case of logistic map), or due to the phenomena found in this paper, this task becomes very difficult if not impossible.
\end{itemize}

\section{Conclusion}
In this paper the fractional order logistic map in the Caputo sense is presented. Even without an analytical proof, it is shown numerically that the iterations number necessary to unveil the real behavior of the FOLM should be of the order of thousands and not of hundreds. Also, it is shown numerically that the bifurcation diagram of the FOLM in the space $(p,x)$ or $(q,x)$ is not unique: for every initial condition another BD is obtained, regardless the iterations number. Therefore, the following natural question arises: which of these bifurcation diagrams should be considered (for theoretical or numerical approach)? Since it seems that this phenomenon appears in other discrete Caputo-like FO systems too, and there is no analytical proof yet to sustain the numerical results in this paper, the subject deserves further investigations, both numerically and even theoretically. To note that a similar phenomenon, related to the integration step-size, has been discovered also in continuous systems of FO \cite{jj}.


\textbf{Funding }No funding to declare.

\textbf{Declaration of Competing Interest} Author declare that he does has no conflict of interest.

\newpage
 \section*{APPENDIX}

\lstset{style=Matlab-editor}
\begin{lstlisting}
function [u]=FO_discrete(q,p,u0,n_max)
%
% Matlab program for solution of the Caputo-like 
% discrete fractional difference Initial Value Problem
%
% \Delta_*^q u(t)=f(t+q-1),u(t+q-1)), t\in N_{1-q}, u(0)=u_0
% 
% f continuous function depending on a real parameter p
% INPUT:
% q -fractional order, q\in(0,1);
% p -parameter;
% u0 -initial condition;
% n_max -maximal number of itereations
% OUTPUT:
% u -solution
% Author: Marius-F. Danca 2021
% web:http://danca.rist.ro/
% email: danca@rist.ro
    u=zeros(n_max,1); % memory alocation
    u(1)=u0+f(u0); % u(1) outside of the loops 
    for n=2:n_max; % $u(2),u(3),...,n(n_max)$
        s=exp(gammaln(n-1+q)-gammaln(n))*f(u0); 
        for j=2:n;
            s=s+exp(gammaln(n-j+q)-gammaln(n-j+1))*f(u(j-1));
        end
        u(n)=u0+s/gamma(q);
    end
    function ff=f(x); % logistic map
    ff=p*x*(1-x);
    end
end
\end{lstlisting}

\newpage{\pagestyle{empty}\cleardoublepage}

\newpage{\pagestyle{empty}\cleardoublepage}

\pagebreak

\end{document}